\documentclass[12pt]{article}
\usepackage{epsfig}

\author{Hugh Nelson Howards \\ Wake Forest University \\
Department of Mathematics \\ Winston Salem, NC 27109 \\
e-mail: howards@wfu.edu}

\title{Limits Of Incompressible Surfaces}

\begin{document}
\maketitle

\newtheorem{definition}{Definition}[section]
\newtheorem{corollary}[definition]{Corollary}
\newtheorem{corrolary}[definition]{Corollary}
\newtheorem{theorem}[definition]{Theorem}
\newtheorem{lemma}[definition]{Lemma}
\newtheorem{claim}[definition]{Claim}
\newtheorem{ex}[definition]{Example}
\newtheorem{q}[definition]{Question}
\newtheorem{exer}[definition]{Exercise}

\begin{abstract}
One can embed arbitrarily many disjoint, non-parallel, non-boundary parallel,
incompressible surfaces in any three manifold
with at least one boundary component of genus two or greater [4].
This paper proves the contrasting, but not contradictory result
that although one can sometimes
embed arbitrarily many surfaces in a 3-manifold it is impossible to
ever embed an infinite number of such surfaces in any compact,
orientable 3-manifold $M$.

\medskip


\medskip
 

\end{abstract}


\section{Introduction and Definitions}

The author would like to thank Michael Freedman and Ying-Qing Wu for 
their helpful comments and advice.

We begin by reviewing a few definitions
 which can be found in most introductory
texts on 3-manifolds.  We rely heavily on Hempel's 
versions in [1].

A surface $(F, \partial{F})$ that is embedded
in a 3-manifold $(M, \partial{M})$ is \underline{properly} 
\underline{embedded} if 
$F \bigcap \partial M = \partial F$.
From this point on when we refer to a surface in a three-manifold,
we will be talking about properly embedded surfaces unless otherwise
noted.
Two surfaces $F_1$ and $F_2$ in a 3-manifold $M$ are
\underline{parallel} if they co-bound a product 
($F \times I; F \times 0 = F_1, F \times 1 = F_2$) in $M$
and $\partial F \times I \subset \partial M$.
A surface $F_1$ in a 3-manifold $M$ is
\underline{boundary parallel} if it co-bounds a product with 
$F_2$, a subsurface of the boundary
($F \times I; F \times 0 = F_1, F \times 1 = F_2$)
and $\partial F \times I \subset \partial M$ in $M$.
A surface $F$ embedded in a three-manifold $M$ is called
\underline{compressible} if any of the following apply.

\begin{enumerate}

\item $F$ is a 2-sphere which bounds a ball in $M$,

\item  $F$ is a disk and either $F \subset \partial M$ or
there is a ball $B \subset M$ such that $\partial B =
F \bigcup D$ where $D$ is a disk contained in $\partial M$, or

\item there is a disk $D \subset M$ with $D \bigcap F = \partial D$
and $\partial D$ not contractible in $F$.  (Note that $D$, of course,
is not required to be properly embedded in $M$).

\end{enumerate}

Otherwise, $F$ is \underline{incompressible}.

A surface $F$ is \underline{boundary compressible} in 
a three-manifold $M$ if either

\begin{enumerate}

\item  $F$ is parallel to a disk in the boundary of $M$ or

\item  There exists a disk $D$ in $M$
such that $D \cap F = \alpha$, an arc in $\partial D$, and
$D \cap \partial M = \beta$ is an arc in $\partial D$ with
$\alpha \cap \beta = \partial \alpha = \partial \beta$ and
$\alpha \cup \beta = \partial D$, and either $\alpha$ ($\beta$) does 
not separate $F$ ($\partial M - \partial F$)
 or $\alpha$ separates F ($\partial M - \partial F$ into two components
and the closure of neither is a disk. 

 Otherwise, F is \underline{boundary incompressible.}

 A 3-manifold, $M$ is 
\underline{irreducible} if every embedded 2-sphere
in $M$ bounds a 3-ball. 
\end{enumerate}  We end with an algebraic definition.
Let $G = G_1 *_H G_2$ designate the \underline{free product}
\underline{with amalgamation}
of the groups $G_1$ and $G_2$ over the group $H$.

\section{The Free Product with Amalgamation}

[4]  demonstrates
that the free group on two generators may be split into
a free product with amalgamation over two arbitrarily large free groups.
This section proves the following contrasting 
result (it is probably known but does not seem to 
have ever been written down):


\begin{theorem}
Let $G_1$ be a group that is not finitely generated  
and let $H$ be finitely generated, then
 $G = G_1 *_H G_2$ is not finitely generated.
\label{fpwa2}
\end{theorem}

\emph{Proof:} \hspace{.5mm}Let $M_i$ be a $K(\pi , 1)$ for $G_i$
$(i=0$ or $1)$.  Let $S_i$ be the image of
$S$ a complex that maps in to represent $H$ in each $M_i$.
Connect $S_0$ to $S_1$ with a cylinder $A = S \times I$ 
with $S \times 0 = S_0$ and $S \times 1 = S_1$ yielding 
a new space, $M$.  
This gives us the free product with
amalgamation for which we are searching.

Choose a base point on $S' = S \times 1/2$
 for $\pi_1 (M)$.  If  $\pi_1 (M)$
is finitely generated, then choose a set of generators,
$\{ \alpha_1, \alpha_2, \dots, \alpha_{i-1}, \alpha_{i}, \dots
\alpha_n \}$,
where $\{ \alpha_1, \alpha_2, \dots, \alpha_{i-1} \}$
generate $H$.
Choose generators for $\pi_1{(M_1)}$
$\{ \gamma_1, \dots \gamma_{i-1}, \gamma_{i}, \dots \}$
and generators $\pi_1{(M_2)}$, 
$\{ \beta_1, \dots \beta_{i-1}, \beta_{i}, \dots \}$
where $\{ \gamma_1, \dots \gamma_{i-1} \}$
(and $\{ \beta_1, \dots \beta_{i-1}\}$)   are 
$\{ \alpha_1, \alpha_2, \dots, \alpha_{i-1} \}$
pushed off into $M_1$ (and $M_2$ respectively).


Now, $\{ \alpha_{i}, \dots \alpha_n \}$
may weave back and forth between $M_1$ and 
$M_2$ a finite number of times and may be expressed as
$$\{ (\gamma_{i_1}^{\pm}
 \beta_{i_1}^{\pm} \dots \gamma_{i_j}^{\pm} \beta_{i_j}^{\pm}),
\dots,(\gamma_{n_1}^{\pm} \beta_{n_1}^{\pm} \dots \gamma_{n_k}^{\pm} \beta_{n_k}^{\pm}) \}$$
where each $\gamma_{s_m}, i \leq s \leq n$ (or $\beta_{s_m}$)
is either some $\gamma_{i}$ or $e$
(or some $\beta_{i}$ or $e$)
Examine any $\gamma_l$ in $\pi_1(M_1)$.
It must be in the span of $\{ \gamma_1, \dots
,\gamma_{i-1}, \gamma_{i_1}, \dots ,\gamma_{n_k} \}$.

To see this, take a 
disk $D$ bounded by  $\gamma_l$ followed by  $\gamma_l^{-1}$
expressed as a product of the generators of $\pi_1(M)$
 that intersects
$S'$ transversally
and minimally.  Since the fundamental group of $S'$
injects, we can assume $ D \cap S'$ has no simple closed curves.
Examine the portion of $ D - D \cap S'$ that contains  $\gamma_l$.
This must be a disk whose boundary consists of generators of 
 $\pi_1(M_1)$ and curves on $S'$ that hence may 
be expressed as products of  $\{ \gamma_1, \dots \gamma_{i-1} \}$ and their
inverses
Therefore, $\gamma_l$ may be written in terms of $\{ \gamma_1, \dots
,\gamma_{i-1}, \gamma_{i_1}, \dots ,\gamma_{n_k} \}$ and their inverses.
Thus, $\pi_1(M_1)$ must be finitely generated.

\begin{corollary}
\label{fpwa3}
Given a compact orientable 3-manifold $M$ and a set of surfaces
$\{ F_i \}$ in $M$, each of the regions in $M - \cup \{ F_i \}$
has finitely generated fundamental group.
\end{corollary}
\emph{Proof:} \hspace{.5mm}After splitting $M$ along a finite number
of the surfaces we attain $M'$ a compact 3-manifold with 
finitely generated fundamental group for which each of the remaining
$\{ F_i \}$ is separating.  None of the complementary regions could 
have infinitely generated fundamental group or else $M'$
would have to, also.


\begin{lemma}
\label{lemma:boundarygenus}
The boundary of a 3-manifold $M$ with finitely generated $\pi_1 (M)$
must have bounded genus.
\end{lemma}
{\em Proof:} \hspace{.5mm} Take a Scott core $C$ for $M$ and expand
$C$ to remain a compact submanifold of $M$, 
but to include an arbitrarily large portion
of the boundary (for example, one might take the $C'$ equal to $C$
plus the closure of a neighborhood of $B \cup A$
where $B$ is a (topologically) large portion
of the boundary and $A$ consists of arcs running from $B$ to $C$).
$C'$ can be further expanded to become a Scott core $C''$ by adding 
2-handles that kill any added elements of $\pi_1$.  Now since
$C''$ has the same fundamental group as $M$ and $H_1$ is just 
abelianized $\pi_1$, $H_1 (C'')$ must have no more generators than 
$\pi_1 (M)$, but the boundary of $C''$ has arbitrarily high genus
($\partial (M) \cap C'' = \partial (M) \cap C'$) and since
$C''$ is compact, ``half lives- half dies'' assures us that 
$H_1 (C'')$ has arbitrarily many generators, which is a contradiction.


\emph{Note:} \hspace{.5mm}  This proof actually shows that any compact
submanifold of $M$ is contained in a Scott core.  We should also point out
that one can also use a 
more complicated homology argument to prove the lemma.

\section{The Behavior of Incompressible Surfaces in a 3-Manifold}

For a while it was claimed that
one could
never embed arbitrarily many disjoint, non-parallel, non-boundary
parallel, incompressible surfaces in 
a three manifold.  The first counter
example was found in [5].  More recently a more general argument
has been used to show that 
any manifold with at least one boundary
component of at least genus two allows such embeddings [4].  
On the other hand, Benedict and
Mike Freedman showed that in any
manifold, if the Euler Characteristic of the
surfaces is bounded, then the number of surfaces will also be bounded [3].
The result in [4]
contrasts with, but does not contradict the main theorem of this paper.

\begin{theorem}
No compact, orientable
three manifold can support an infinite number of disjoint, non-parallel,
non-boundary-parallel, incompressible surfaces.
\label{infsurf}
\end{theorem}

{\em Note:} \hspace{.5mm}An easy argument 
in [3] shows that any given three manifold supports
only a bounded number of boundary-parallel, but non-parallel surfaces, if
it is assumed that none of the surfaces are disks or annuli. Thus, this
assumption could replace the non-boundary parallel assumption above.

{\em Proof:} \hspace{.5mm} To begin our proof,
 we should recall that the usual Haken
finiteness argument using normal
surfaces shows that there can only be a finite number of incompressible,
boundary-incompressible surfaces, so we need only consider the surfaces
that are incompressible but boundary compressible. 

Given this, in order to derive a contradiction we may assume
that we have an infinite list of surfaces $\{F_i\}$.
It will be convenient later to assume we have no annuli,
and [3] assures us that $M$ can only have a finite number of possible disjoint
non-boundary-parallel annuli,
therefore that we may assume that none of the
surfaces on the infinite list of $\{F_i\}$
are annuli.
Note that $M$ only has a finite
number of
boundary components.  Also note that since each $F_i$ is boundary
compressible,
each one meets at least one boundary component of $M$ in a set of simple
closed
curves.  We will choose one of the surfaces and examine its boundary
compression disk.

We examine how the $\{F_i\}$ intersect the boundary of $M$ and
define a product region in $\partial M$ to be an annular
region with boundary two parallel curves from the boundary of
the $\{F_i\}$.  Of course a non-product region simply refers to
a region of the boundary which is not a product region.  See
Figure 1.

\begin{figure}[htbp] 
\centerline{\hbox{\epsfig{file=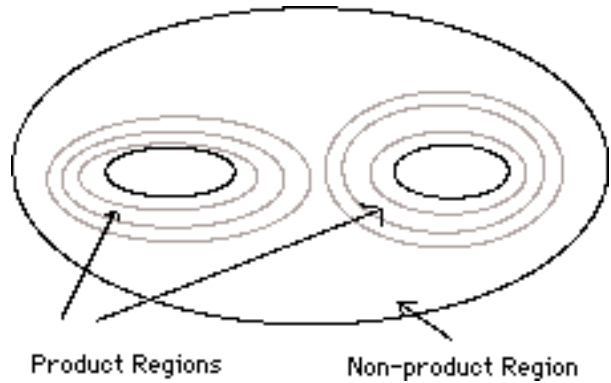}}}
\caption{Product and non-product regions}\label{fig:prodnonprod}
\end{figure}

\begin{lemma}
The boundary components of the surfaces $F_i$ are each parallel
to one of a finite number of curves on the boundary of $M$.
\end{lemma}
{\em Proof:} \hspace{.5mm}  Since $M$ has
 a finite number of boundary components and
all of the $F_i$ are disjoint and therefore have disjoint boundary
components, the proof is an easy
application of hierarchy arguments to non-trivial
simple closed curves on closed orientable surfaces. 

\vspace{.25in}

Since there are only a finite number of curve types on the boundary
there can only be a finite number of non-product regions.  These in kind
can only correspond to $\{M_1, \dots, M_n\}$, a finite subset
of the regions obtained by cutting $M$
up along the union of the $\{F_i\}$.

Let $M'_i$ be a closed regular neighborhood in $M_i$ of the boundary
of $M_i$.  This means that $M'_i$ is a (not necessarily
connected) surface crossed with the unit interval and therefore has
incompressible boundary.


Now in $M$ replace $M_i$ by $M'_i$ obtaining $M'$.    
Since we have just observed that 
the $\{ M_i \}$ is a finite set of pieces and the pieces that do not have
non-product regions on their boundary are left alone and not replaced,
we have only altered a finite number of regions.
Corollary~\ref{fpwa3} together with Lemma~\ref{lemma:boundarygenus}
assures us that for any region $M_i, \partial M_i$ has bounded genus.
We also note that if the boundary of a 3-manifold has a ``puncture,"
the puncture has to extend to infinity, so a neighborhood of the puncture
must be an infinitely long annulus.  Since the boundary of
$M_i$ is made from a (potentially infinite) list of compact surfaces,
such a puncture could only result from an infinite list of annuli
glued together, but this is impossible since 
none of the $\{F_i\}$ are annuli.
$M$ was compact, and each of the new pieces are compact, 
therefore $M'$ is, too.


The infinite
collection of disjoint, non-parallel,
non-boundary-parallel, incompressible surfaces in $M$ become
an infinite
collection of disjoint, non-parallel,
non-boundary-parallel, incompressible surfaces in $M'$.
The remainder of the paper will show that the surfaces are 
also boundary-incompressible in $M'$ which is a contradiction.

We now choose an $F_i$ and examine a disk representing its
boundary compression.  Choose the outermost such 
boundary compression disk with respect to $F_i$, so that the interior
of the disk is disjoint from $F_i$.
We look at the intersection of the
disk with the other surfaces.  Since the surfaces are incompressible,
we may use an innermost loop argument to show that we can choose
to have only arcs and no simple closed curves in the intersection.


\begin{lemma}
\label{nonprod1}
The boundary of boundary compressing disks runs through a
non-product region.
\end{lemma} 
{\em Proof:} \hspace{.5mm} If the
 component were strictly in product regions, either
the component would connect a boundary component 
to itself in a trivial manner
which is prohibited by the 
definition of a boundary compression, or it would cross
the annulus in the unique arc connecting the two boundary components.
(See Figure 2.)

\begin{figure}[htbp] 
\centerline{\hbox{\epsfig{figure=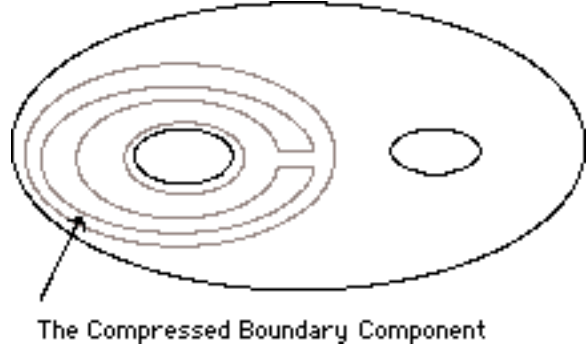}}}
\caption{A boundary compression in a product region}\label{fig:bdrycompr}
\end{figure}

In the latter case, after compressing
we have a boundary component which is trivial in the fundamental
group of $\partial M$.  So now we have two options:
either we have a disk, which is impossible since that
would mean we started with an annulus,
or else we have a compressible surface.  This is also a contradiction
as performing a boundary compression cannot make a surface compressible that
was not already compressible.  


 The boundary compressing disk must therefore intersect
a non-product region essentially.  
This yields a compressing disk
for a boundary component of one of the $M'_i$
(since those are the only pieces which can have a non-product
region on its boundary), but these pieces have incompressible 
boundary, so this is a contradiction.  There is
no boundary compressing disk and therefore the surfaces must be boundary
incompressible in $M'$.

Haken finiteness now applies, so there are only a finite
number of surfaces. This contradicts our original
assumption.

\pagebreak
\section{References}

\hspace{.2in}

[1]  Hempel, J., \it{ 3-Manifolds} 
(Annals of Mathematics Studies, 86, Princeton
University Press, Princeton, NJ, 1976).

\hspace{.2in}

[2]  Jaco, W., \it{ Lecture Notes on 3-manifold Topology}
(CBMS Regional Conf. Ser. in Math.,
No. 43 Amer. Math. Soc., Providence, RI, 1980).

\vspace{.25in}

[3]  Freedman, B.; Freedman, M.. Kneser-Haken finiteness for 
bounded $3$-manifolds locally free groups, and cyclic covers. 
Topology 37 (1998), no. 1, 133--147. 

\vspace{.25in}

[4]  Howards, H,  Arbitrarily many incompressible surfaces
in three-manifolds with genus two boundary components: preprint 1996

\vspace{.25in}

[5]  Sherman, W,  Ph.D Thesis UCLA

\end{document}